\documentclass[12pt]{article}
\usepackage{e-jc}

\usepackage[english]{babel}
\usepackage{amsthm,amsmath,amssymb}
\usepackage{graphicx}
\usepackage[small,bf,hang]{caption} 

\usepackage{algorithm}                   
\usepackage{algorithmic}

\usepackage{epsf} 

\usepackage{hyperref}
\usepackage{enumerate} 

\usepackage{multirow}

\newtheorem{theorem}{Theorem}
\theoremstyle{plain}

\newtheorem{proposition}[theorem]{Proposition}

\theoremstyle{definition}

\graphicspath{{figures/}}

\title{\bf On the smallest snarks with \\oddness 4 and connectivity 2}

\author{Jan Goedgebeur\thanks{Supported by a Postdoctoral Fellowship of the Research Foundation Flanders (FWO).}\\
\small Department of Applied Mathematics, Computer Science and Statistics\\[-0.8ex]
\small Ghent University\\[-0.8ex]
\small Krijgslaan 281-S9,\\[-0.8ex]
\small 9000 Ghent, Belgium\\[-0.8ex]
\small\tt jan.goedgebeur@ugent.be\\
\\
}

\date{\dateline{XX}{XX}}

\begin{document}

\maketitle

\begin{abstract}
A \textit{snark} is a bridgeless cubic graph which is not 3-edge-colourable. The \textit{oddness} of a bridgeless cubic graph is the minimum number of odd components in any 2-factor of the graph.

Lukot'ka, M{\'a}{\v{c}}ajov{\'a}, Maz{\'a}k and \v{S}koviera showed in [\emph{Electron.\ J.\ Combin.} 22 (2015)] that the smallest snark with oddness 4 has 28 vertices and remarked that there are exactly two such graphs of that order. 
However, this remark is incorrect as -- using an exhaustive computer search -- we show that there are in fact three snarks with oddness 4 on 28 vertices. In this note we present the missing snark and also determine all snarks with oddness 4 up to 34 vertices.

\end{abstract}


\section{Introduction and main result}

The \textit{chromatic index} of a graph $G$ is the minimum number of colours required for an edge colouring of that graph such that no two adjacent edges have the same colour. It follows from Vizing's classical theorem that a cubic graph has chromatic index either 3 or 4.
Isaacs~\cite{isaacs_75} called cubic graphs with chromatic index 3 \textit{colourable} and those with chromatic index 4 \textit{uncolourable}.
Cubic graphs with bridges can easily be seen to be uncolourable and are therefore considered to be trivially uncolourable.

A \textit{snark} is a bridgeless cubic graph which is not 3-edge-colourable. Note that in the literature stronger conditions are sometimes required for a graph to be a snark, e.g.\ that it also must have girth at least 5 and be cyclically 4-edge-connected. (The \textit{girth} of a graph is the length of its shortest cycle and a graph is \textit{cyclically $k$-edge-connected} if the deletion of fewer than $k$ edges from the graph does not create two components both of which contain at least one cycle). Here we will focus on snarks with girth at least 4 since snarks with triangles can be easily reduced to smaller triangle-free snarks. 

One of the reasons why snarks are interesting is the fact that the smallest counterexamples to several important conjectures (such as the cycle double cover conjecture~\cite{Sey79,Sze73} and the 5-flow conjecture~\cite{tutte_5-flow}) would be snarks.

The \textit{oddness} of a bridgeless cubic graph is the minimum number of odd components in any 2-factor of the graph. The oddness is a natural measure for how far a graph is from being 3-edge-colourable. It is straightforward to see that a bridgeless cubic graph is 3-edge-colourable if and only if it has oddness 0. Also note that the oddness must be even as cubic graphs have an even number of vertices.

Snarks with large oddness are of special interest since several conjectures (including the cycle double cover conjecture and the 5-flow conjecture) are proven to be true for snarks with small oddness.

In~\cite{lukotka2015small} Lukot'ka, M{\'a}{\v{c}}ajov{\'a}, Maz{\'a}k and \v{S}koviera showed the following.

\begin{theorem}[Theorem~12 in~\cite{lukotka2015small}]\label{thm:oddness4}
The smallest snark with oddness 4 has 28 vertices. There is one such snark with cyclic connectivity 2 and one with cyclic connectivity 3.
\end{theorem}

After the proof of this theorem they remark:

\begin{quote}
``The computer searches referred to in the proof of Theorem~12 can be extended to show that there are exactly two snarks with oddness 4 on 28 vertices -- those displayed in Figure 2 (from~\cite{lukotka2015small})." 
\end{quote}

However, this remark is incorrect as -- using an exhaustive computer search -- we have shown that there are in fact three snarks with oddness 4 on 28 vertices, which leads to the following proposition.

\begin{proposition}
There are exactly three snarks with oddness 4 on 28 vertices. There are two such snarks with cyclic connectivity 2 and one with cyclic connectivity 3.
\end{proposition}

\begin{figure}[h!tb]
	\centering
	\includegraphics[width=0.8\textwidth]{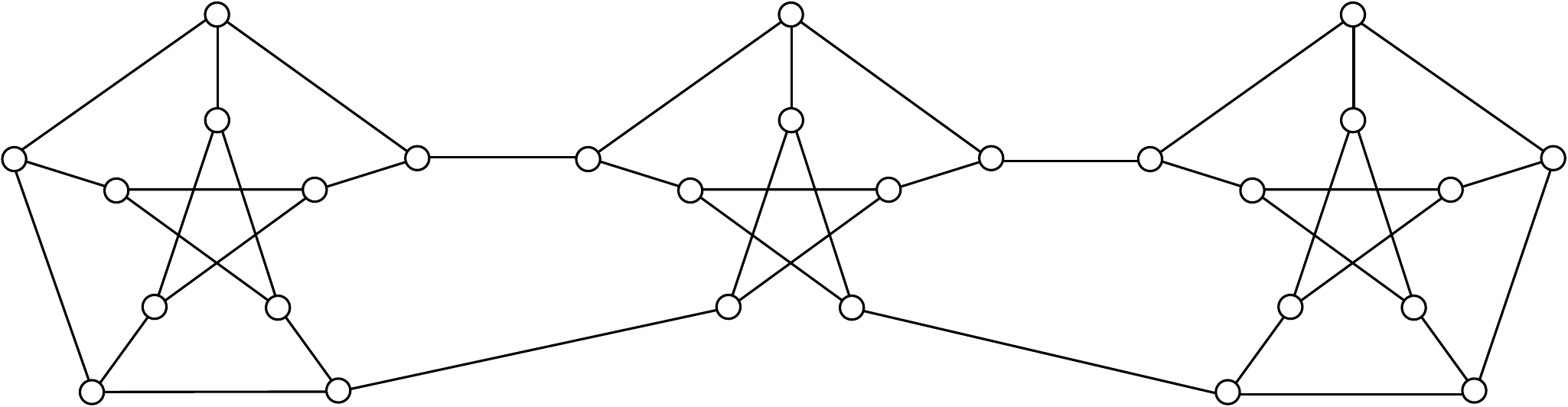}
	\caption{The snark with oddness 4 on 28 vertices which was missing in~\cite{lukotka2015small}.}
	\label{fig:28_missing}
\end{figure}

The missing snark from Theorem~\ref{thm:oddness4} has connectivity 2 and is shown in Figure~\ref{fig:28_missing}.
Using the program \textit{snarkhunter}~\cite{snark-paper,brinkmann_11} we have generated all snarks with girth at least 4 up to 34 vertices and tested which of them have oddness greater than 2.  
The results of this search are listed in Table~\ref{table:counts_girth4} and in Table~\ref{table:counts_girth5} for snarks with girth at least 5. In~\cite[Lemma~2]{lukotka2015small} it was shown that the smallest snarks of a given oddness have girth at least 5.

\begin{table}
\centering
\setlength{\tabcolsep}{8pt} 
 \renewcommand{\arraystretch}{1.15} 
	\begin{tabular}{c r ccc }
		 \multirow{2}{*}{Order} & \multirow{2}{*}{All} & \multicolumn{3}{c}{Oddness 4}\\
		\cline{3-5}
		  & & Connectivity 2 & Connectivity 3 & Total  \\
		 \hline
10 & 1 & 0 & 0 & 0 \\
12 & 0 & 0 & 0 & 0 \\
14 & 1 & 0 & 0 & 0 \\
16 & 4 & 0 & 0 & 0 \\
18 & 26 & 0 & 0 & 0 \\
20 & 167 & 0 & 0 & 0 \\
22 & 1 448 & 0 & 0 & 0 \\
24 & 15 168 & 0 & 0 & 0 \\
26 & 189 861 & 0 & 0 & 0 \\
28 & 2 716 555 & 2 & 1 & 3 \\
30 & 43 504 872 & 9 & 4 & 13 \\
32 & 767 442 160 & 57 & 32 & 89 \\
34 & 14 752 529 374 & 454 & 313 & 767 \\		
	\end{tabular}

\caption{The counts of all 2-connected snarks with girth at least 4 up to 34 vertices and the number of snarks with oddness 4 among them.}
\label{table:counts_girth4}
\end{table}

\begin{table}
\centering
\setlength{\tabcolsep}{8pt} 
 \renewcommand{\arraystretch}{1.15} 
	\begin{tabular}{c r ccc }
		 \multirow{2}{*}{Order} & \multirow{2}{*}{All} & \multicolumn{3}{c}{Oddness 4}\\
		\cline{3-5}
		  & & Connectivity 2 & Connectivity 3 & Total  \\
		 \hline
10 & 1 & 0 & 0 & 0 \\
12 & 0 & 0 & 0 & 0 \\
14 & 0 & 0 & 0 & 0 \\
16 & 0 & 0 & 0 & 0 \\
18 & 3 & 0 & 0 & 0 \\
20 & 14 & 0 & 0 & 0 \\
22 & 107 & 0 & 0 & 0 \\
24 & 1 109 & 0 & 0 & 0 \\
26 & 15 255 & 0 & 0 & 0 \\
28 & 236 966 & 2 & 1 & 3 \\
30 & 4 043 956 & 9 & 4 & 13 \\
32 & 74 989 646 & 33 & 21 & 54 \\
34 & 1 500 084 086 & 139 & 138 & 277 \\
	\end{tabular}

\caption{The counts of all 2-connected snarks with girth at least 5 up to 34 vertices and the number of snarks with oddness 4 among them.}
\label{table:counts_girth5}
\end{table}

None of the snarks up to 34 vertices has oddness greater than 4, so this yields the following proposition.

\begin{proposition}
The smallest 2-connected snark with oddness at least 6 has at least 36 vertices.
\end{proposition}

All snarks with oddness at least 4 up to at least 36 vertices have (cyclic) connectivity 2 or 3, since it follows from~\cite{snark-paper} and~\cite{GMS} that there are no cyclically 4-edge-connected snarks with oddness at least 4 up to at least 36 vertices.

We implemented two independent algorithms to compute the oddness of a bridgeless cubic graph: the first algorithm computes the oddness by constructing perfect matchings while the second algorithm does this by constructing 2-factors directly by searching for disjoint cycles. The source code of both programs can be obtained from~\cite{oddness-site}. All of our results reported in this article were independently confirmed by both programs. 

The graphs from Tables~\ref{table:counts_girth4} and~\ref{table:counts_girth5} can be downloaded and inspected in the database of interesting graphs from the \textit{House of Graphs}~\cite{hog} by searching for the keywords ``snark * oddness 4''.

\begin{figure}[h!tb]
	\centering
	\includegraphics[width=0.9\textwidth]{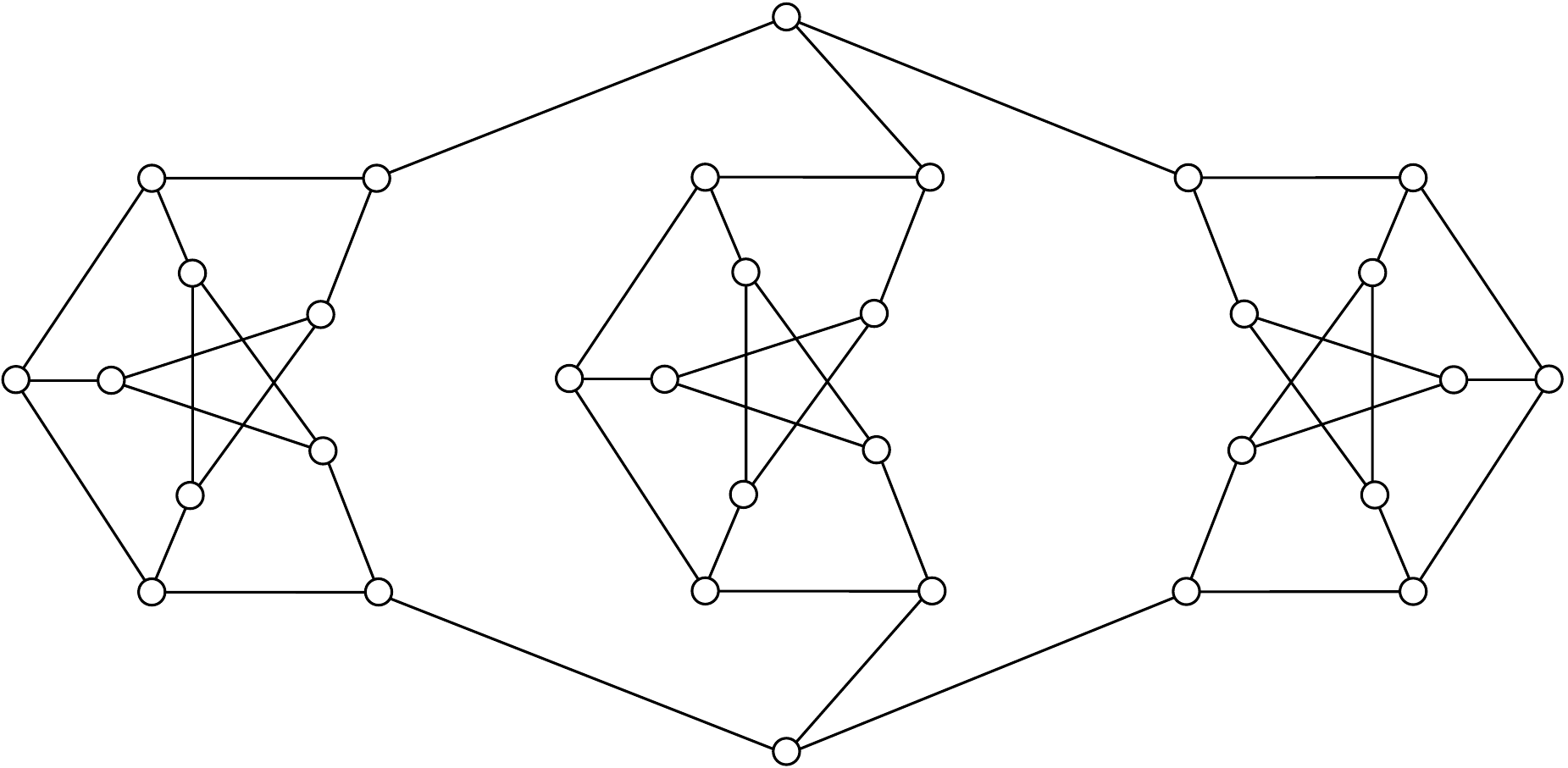}
	\caption{The most symmetric snark with girth at least 4 and with oddness 4 up to 34 vertices. It has 32 vertices and its automorphism group has order~768.}
	\label{fig:32v_most_symm}
\end{figure}

The most symmetric snark with girth at least 4 and with oddness 4 up to 34 vertices is shown in Figure~\ref{fig:32v_most_symm}. It has 32 vertices, connectivity 2, girth 5 and its automorphism group has order~768.

\subsection*{Acknowledgements}

We would like to thank Martin \v{S}koviera for useful suggestions.
Most of the computations were carried out using the Stevin Supercomputer Infrastructure at Ghent University.



\begin{thebibliography}{1}

\bibitem{hog}
G.~Brinkmann, K.~Coolsaet, J.~Goedgebeur, and H.~M{\'e}lot.
\newblock {House of Graphs: a database of interesting graphs}.
\newblock {\em Discrete Applied Mathematics}, 161(1-2):311--314, 2013.
\newblock Available at \url{http://hog.grinvin.org/}.

\bibitem{snark-paper}
G.~Brinkmann, J.~Goedgebeur, J.~H{\"a}gglund, and K.~Markstr{\"o}m.
\newblock Generation and properties of snarks.
\newblock {\em Journal of Combinatorial Theory, Series B}, 103(4):468--488,
  2013.

\bibitem{brinkmann_11}
G.~Brinkmann, J.~Goedgebeur, and B.~D. McKay.
\newblock Generation of cubic graphs.
\newblock {\em Discrete Mathematics and Theoretical Computer Science},
  13(2):69--80, 2011.

\bibitem{oddness-site}
J.~Goedgebeur.
\newblock Source code of two programs to compute the oddness of a graph:
  \url{http://caagt.ugent.be/oddness/}.
  
\bibitem{GMS} J.~Goedgebeur, E.~M{\'a}{\v{c}}ajov{\'a}, and M.~{\v{S}}koviera.
 \newblock  Smallest snarks with oddness 4 and cyclic connectivity 4 have order 44.
 \newblock arXiv:1712.07867, 2017.  

\bibitem{isaacs_75}
R.~Isaacs.
\newblock {Infinite families of nontrivial trivalent graphs which are not Tait
  colorable}.
\newblock {\em American Mathematical Monthly}, 82(3):221--239, 1975.

\bibitem{lukotka2015small}
R.~Lukot'ka, E.~M{\'a}{\v{c}}ajov{\'a}, J.~Maz{\'a}k, and M.~{\v{S}}koviera.
\newblock Small snarks with large oddness.
\newblock {\em The Electronic Journal of Combinatorics}, 22(1), 2015.

\bibitem{Sey79}
P.~D. Seymour.
\newblock Sums of circuits.
\newblock In {\em Graph theory and related topics ({P}roc. {C}onf., {U}niv.
  {W}aterloo, {W}aterloo, {O}nt., 1977)}, pages 341--355. Academic Press, New
  York, 1979.

\bibitem{Sze73}
G.~Szekeres.
\newblock Polyhedral decompositions of cubic graphs.
\newblock {\em Bulletin of the Australian Mathematical Society}, 8:367--387,
  1973.

\bibitem{tutte_5-flow}
W.~T. Tutte.
\newblock A contribution to the theory of chromatic polynomials.
\newblock {\em Canadian Journal of Mathematics}, 6:80--91, 1954.

\end{thebibliography}
\end{document}